\documentclass[preprint,3p,11pt]{elsarticle}

\usepackage{tabls}
\usepackage{cites}
\usepackage{epsf}
\usepackage{appendix}
\usepackage{ragged2e}

\usepackage[ruled,vlined]{algorithm2e}
\usepackage[caption = false]{subfig}
\usepackage{amsmath}
\usepackage{amssymb}
\usepackage{xcolor}

\usepackage{hyperref}
\hypersetup{colorlinks=true}

\newcommand{\dv}[2]{\frac{\partial #1}{\partial #2}}
\newcommand{\ip}[1]{\langle #1 \rangle}

\begin{document}
\begin{frontmatter}



\title{Multilevel  Method for Thermal Radiative Transfer Problems
	with Method of Long Characteristics for the  Boltzmann  Transport  Equation}
\setcounter{ead}{2}
\author[lanl]{Joseph M. Coale}
\ead{jmcoale@lanl.gov}
\author[ncsu]{Dmitriy Y. Anistratov}
\ead{anistratov@ncsu.edu}
\address[lanl]{Los Alamos National Laboratory,
  Los Alamos, NM  87545}
\address[ncsu]{Department of Nuclear Engineering,
North Carolina State University,
 Raleigh, NC 27695}

\begin{abstract}
  In this paper analysis is performed on a computational method for thermal radiative transfer (TRT) problems based on the multilevel
  quasidiffusion (variable Eddington factor) method with
  the method of long characteristics
  (ray tracing)
  for the Boltzmann transport equation (BTE).
  The method is formulated with a multilevel set of moment equations of the BTE which are coupled to the material energy balance (MEB).
  The moment equations are exactly closed via the Eddington tensor defined by the  BTE solution.
  Two discrete spatial meshes are defined: a material grid on which the MEB and low-order moment equations are
  discretized, and a grid of characteristics for solving
  the BTE.
  Numerical testing of the method is completed on the well-known Fleck-Cummings test problem which models a supersonic radiation wave propagation.
  Mesh refinement studies are performed on each of the two spatial grids independently, holding one mesh width constant while refining the other.
  We also present the data on convergence of iterations.
\end{abstract}


 \end{frontmatter}

\section{Introduction}
In this paper we investigate the performance of computational methods for high-energy density physics (HEDP)
simulations \cite{zel-1966,drake-2006}. Phenomena that occur in the high-energy density (HED) regime are typically modeled with complex multiphysical systems of partial differential equations (PDEs).
Energy redistribution in the HED regime is dominated by the mechanism of radiative transfer, where energy is transported by the emission and absorption of photon radiation.
The radiative transfer equation is the Boltzmann transport equation (BTE)
which models how photons propagate through and interact with matter.
Thus, the BTE is
an essential component in the models of HED phenomena.
These systems of PDEs
exhibit several properties which present challenges for their numerical simulation, including: (i) strong nonlinearity, (ii) tight coupling between equations, (iii) multiscale behavior in space-time, (iv) high-dimensionality.

Discretization of the equations involved in HEDP simulations must be handled with these properties in mind to preserve the fundamental physical behavior of the modeled phenomena.
To this end
the method of long characteristics (MOLC), a.k.a. ray-tracing schemes (RTS), for the BTE
hold some advantage compared to other discretization schemes.
Numerical techniques based on MOLC/RTS for solving particle transport problems are
formulated by means of
the integral form of the BTE. They belong to the family of the method of characteristics that has a rich and long history of development  \cite{Vladimirov-1958,goldin-1960,takeuchi-1969,askew-1972,brough-chudley-1980,Alcouffe-1981,askew-roth-1982,casmo-1993,kord-2002,Sanchez-2000,zika-adams-2000,Abel-2002,Pandya-2009,Buntemeyer-2016}.

This terminology reflects that the MOLC solves the BTE along characteristics extending over the entire spatial domain. The MOLC doesn’t employ interpolation of solution on faces of cells.
Since the solution is not interpolated on the spatial mesh, discontinuities can be well resolved and carried across the problem domain.
The BTE solution will also be defined on a unique discrete grid separate from the multiphysics equations it becomes coupled to.
It is therefore straightforward to refine or coarsen the discrete BTE grid independently of the multiphysics equations, allowing for computational resources to be more effectively managed.

We consider multi-dimensional, time dependent problems that capture the essential challenges in HEDP simulation.
Specifically, the computational method is developed on and analyzed with the fundamental thermal radiative transfer (TRT) problem, formulated by the multigroup BTE for photons
\begin{equation} \label{eq:bte}
	\frac{1}{c}\dv{I_g}{t} + \boldsymbol{\Omega}\cdot\boldsymbol{\nabla}I_g + \varkappa_g(T)I_g = \varkappa_g(T)B_g(T) \,  ,
\end{equation}
\begin{equation}
	I_g|_{\boldsymbol{r}\in\partial\Gamma}=I_g^\text{in}\ \ \text{for}\ \ \boldsymbol{\Omega}\cdot\boldsymbol{n}_\Gamma<0,\quad I_g|_{t=0}=I_g^0 \, ,
\end{equation}
\begin{equation*}
	\boldsymbol{r}\in\Gamma, \  t> 0, \  \boldsymbol{\Omega}\in\mathcal{S}, \ g=1,\dots,G
\end{equation*}
coupled to the material energy balance (MEB) equation that describes radiation-matter
energy exchange
\begin{equation} \label{meb}
	\dv{\varepsilon(T)}{t} = {\sum_{g=1}^{G}\bigg(\int_{4\pi}I_gd\Omega - B_g(T)\bigg)\varkappa_g(T)},\quad T|_{t=0}=T^0 \, ,
\end{equation}
where $\boldsymbol{r}$ is spatial position, $\boldsymbol{\Omega}$ is the direction of particle motion, $g$ is the photon frequency group index, $c$ is the speed of light, $\Gamma$ is the spatial domain, $\partial\Gamma$ is the domain boundary, $\boldsymbol{n}_\Gamma$ is the unit normal to $\partial\Gamma$ and $\mathcal{S}$ is the unit sphere. $I_g(\boldsymbol{r},\boldsymbol{\Omega},t)$ is the group specific intensity of radiation, $T(\boldsymbol{r},t)$ is the material temperature, $\varkappa_g(\boldsymbol{r},t;T)$ is the group photon opacity, $\varepsilon(\boldsymbol{r},t;T)$ is the material energy density and $B_g(\boldsymbol{r},t;T)$ is the group Planck black-body distribution function.
This TRT problem captures all fundamental challenges and features associated with multiphysical HEDP problems, and serves as a useful platform for the analysis of computational methods aimed at HEDP simulation.

To solve the TRT problem we employ the multilevel quasidiffusion (QD), a.k.a. Variable Eddington Factor (VEF), methodology \cite{gol'din-1964,auer-mihalas-1970,gol'din-1972,winkler-norman-mihalas-85,PASE-1986,mihalas-FRH-1984,dya-aristova-vya-mm1996,aristova-vya-avk-m&c1999,dya-jcp-2019}, which is essentially a nonlinear method of moments.
The multilevel QD (MLQD) method is defined by a system consisting of (i) the high-order BTE and (ii) a hierarchy of low-order equations for moments of the radiation intensity. The low-order equations are exactly closed through the Eddington tensor and other linear-fractional factors that are weakly dependent on the BTE solution.
Multiphysical equations (i.e. the MEB) are coupled to the low-order equations on the scale of the multiphysics, effectively reducing the dimensionality of the problem.
The MOLC/RTS is used to compute the high-order solution of the BTE that is used to compute the Eddington tensor to close the moment equations. The important property of the MLQD method is that the hierarchy of high-order BTE and moment equations
can be discretized independently \cite{dya-vyag-ttsp}. This enables one to take advantage of essentially different type of discretization for high-order and low-order equations to find an improvement in the numerical solution and computational efficiency.

This paper aims  to address
several fundamental questions regarding the use of
MOLC/RTS as
part of the MLQD hierarchy of equations; Specifically as generating the Eddington tensor closures for HEDP applications.
The discretization of the system of high-order and low-order equations are performed on two grids. The low-order equations for moments of the intensity and the MEB equation for temperature are
approximated on the underlying mesh defined by material properties. We refer to it as the spatial material mesh.
The ray-tracing technique in MOLC/RTS defines a mesh of characteristics  over the whole spatial domain.  The characteristics stretch from one domain boundary to another. This procedure also creates a subgrid of
characteristics in each spatial material cell. Note that the  mesh of characteristics and the resulting subgrid in spatial material cells
depend on discrete directions. It can also vary for different frequency groups.
We investigate the essential relation between the two distinct discrete grids by performing studies on the effects of refining each grid independently of one another.
This process is meant to shed light on the relative importance of each individual grid and how much effort should be given to refine either the characteristic mesh or underlying material mesh to optimize solution accuracy and performance.
Furthermore, the impact of mesh sizing on iterative convergence of the solution is examined.

The remainder of the paper is organized as follows.
In Sec. \ref{sec:method}, the  MLQD  method  is formulated.
In Sec. \ref{MOLC}, we present details of a variant of MOLC/RTS for solving the high-order BTE.
The numerical results are presented in Sec. \ref{sec:res}.
We conclude with a discussion in Sec. \ref{sec:conc}.

\section{Multilevel QD/VEF Method for TRT \label{sec:method}}

Two systems of moment equations of the BTE are constructed, the first of which being the multigroup low-order quasidiffusion (LOQD) equations
\begin{subequations} \label{eq:mg-loqd}
	\begin{equation}
		\dv{E_g}{t} + \boldsymbol{\nabla}\cdot\boldsymbol{F}_g + c\varkappa_g(T)E_g = 4\pi\varkappa_g(T)B_g(T) \,  ,
	\end{equation}
	\begin{equation}
		\frac{1}{c}\dv{\boldsymbol{F}_g}{t} + c\boldsymbol{\nabla}\cdot (\boldsymbol{\mathfrak{f}}_g E_g) + \varkappa_g(T) \boldsymbol{F}_g = 0 \, ,
	\end{equation}
\end{subequations}
which solve for the multigroup radiation energy density $E_g=\frac{1}{c}\int_{4\pi}I_gd\Omega$ and flux $\boldsymbol{F}_g=\int_{4\pi}\boldsymbol{\Omega}I_gd\Omega$. The second system of moment equations is the effective grey LOQD equations
\begin{subequations} \label{eq:gr-loqd}
	\begin{equation}
		\dv{E}{t} + \boldsymbol{\nabla}\cdot\boldsymbol{F} + c\ip{\varkappa}_E E = c\ip{\varkappa}_B a_RT^4 \, ,	
	\end{equation}
\begin{equation}
		\frac{1}{c}\dv{\boldsymbol{F}}{t} + c\boldsymbol{\nabla} \cdot (\ip{\boldsymbol{\mathfrak{f}}}_{E}E) +
\bar{\mathbf{K}}_{R}\boldsymbol{F}
 + \bar{\boldsymbol{\eta}}E = 0 \, ,
	\end{equation}
\end{subequations}
which solve for the total radiation energy density $E=\sum_{g=1}^GE_g$ and flux $\boldsymbol{F}=\sum_{g=1}^G\boldsymbol{F}_g$.
The MEB equation is cast in effective grey form
\begin{equation} \label{grey_meb}
	\dv{\varepsilon(T)}{t} = c\ip{\varkappa}_E E - c\ip{\varkappa}_Ba_RT^4,
\end{equation}
and coupled with Eqs. \eqref{eq:gr-loqd}.
The multigroup LOQD equations \eqref{eq:mg-loqd} are exactly closed via the Eddington tensor computed by the solution of the BTE   \eqref{eq:bte}:
\begin{equation} \label{eq:et}
	\boldsymbol{\mathfrak{f}}_g =
	\frac{ \int_{4\pi}\boldsymbol{\Omega}\otimes\boldsymbol{\Omega}{I}_g\ d\Omega}
	{\int_{4\pi}I_gd\Omega } \, .
\end{equation}
The spectrum average opacities and
 coefficients
 calculated by the solution of multigroup LOQD equations \eqref{eq:mg-loqd}
which define exact closures for the LOQD system are
\begin{equation}   \label{eq:grip1}
	\ip{u}_E = \frac{\sum_{g=1}^Gu_gE_g}{\sum_{g=1}^GE_g},\quad
	\bar{\mathbf{K}}_{R}=\text{diag}\big(\ip{u}_{|F_x|},\ip{u}_{|F_y|},\ip{u}_{|F_z|}\big) , ,
\end{equation}
\begin{equation} \label{eq:grip2}
	\ip{u}_{|F_\alpha|} = \frac{\sum_{g=1}^Gu_g|F_{\alpha,g}|}{\sum_{g=1}^G|F_{\alpha,g}|},\quad
	\bar{\boldsymbol{\eta}} = \frac{\sum_{g=1}^{G} (\varkappa_g- 	\bar{\mathbf{K}}_{R})\boldsymbol{F}_g }{\sum_{g=1}^GE_g} \, .
\end{equation}
In sum the MLQD method for the TRT problem is formulated with:
\begin{enumerate}
	\item The BTE \eqref{eq:bte} discretized with the MOLC/RTS,
	\item The multigroup LOQD equations (Eqs. \eqref{eq:mg-loqd} \& \eqref{eq:et}) discretized with the second-order finite volume (FV) scheme \cite{pg-dya-jcp-2020},
	\item The {\it effective grey problem} formed by the closed system of effective grey LOQD equations and the MEB equation (Eqs. \eqref{eq:gr-loqd}, \eqref{grey_meb}, \eqref{eq:grip1} \& \eqref{eq:grip2}). The spatial approximation FV scheme for effective grey LOQD equations are algebraically consistent with the scheme for  the multigroup LOQD equations \cite{pg-dya-jcp-2020}.
\end{enumerate}


\section{Method of Long Characteristics / Ray-Tracing Scheme} \label{MOLC}

The MOLC/RTS is derived by first performing a change of coordinates from $(x,y)$ to $(u,v)$ where $\boldsymbol{e}_u=\boldsymbol{\Omega}$ and $\boldsymbol{e}_u\cdot\boldsymbol{e}_v=0$.
The BTE discretized in time with the backward-Euler scheme, along characteristics is given
\begin{equation}
	\dv{I^n_g(u)}{u} + \tilde{\varkappa}^n_g(u) I_g^n(u) = {Q}^n_g(u) \, ,
	\label{eq:bte_char}
\end{equation}
\begin{equation}
	\tilde{\varkappa}^n_g = \varkappa^n_g + \frac{1}{c\Delta t^n},\quad Q^n_g= \varkappa_g^nB_g^n + \frac{1}{c\Delta t^n}I_g^{n-1} \, ,
\end{equation}
and $\Delta t^n$ is the $n^\text{th}$ time step.
Let all low-order equations be discretized on an `underlying' orthogonal material spatial grid.
The material temperature is a piece-wise function on the set material spatial cells. As a result, the opacities and  Planckian emission source are constant in each  material cell.
With the MOLC/RTS, the BTE is discretized on a characteristic grid traced over the underlying grid for $M$ discrete directions $\boldsymbol{\Omega}_m$.
For each discrete direction, there will be $K_m$ characteristics traced over the spatial domain for a total of $K$ characteristics which formulate the high-order discrete grid.
Let $u^+$ and $u^-$ be the intersection (entry and exit) points of a given characteristic with the boundaries of a single cell in the underlying grid. This forms a segment of the characteristic with length $\ell=(u^+-u^-)\sin(\zeta_m)$ where $\zeta_m$ is the angle between $\boldsymbol{\Omega}_m$ and the $z$ axis. Integrating Eq. \eqref{eq:bte_char} along each characteristic segment gives
\begin{gather}
	I_{k,s+1} = I_{k,s}e^{-\frac{\tilde{\varkappa}_{i}\ell_{k,s}}{\sin(\zeta)}} + \frac{{Q}_{i}}{\tilde{\varkappa}_{i}} \big(1 - e^{-\frac{\tilde{\varkappa}_{i}\ell_{k,s}}{\sin(\zeta)}}\big),\quad
	{\bar{I}}_{k,s} = \alpha_{k,s}I_{k,s} + (1-\alpha_{k,s})I_{k,s+1}\, ,\\
	k=1,\dots,K,\quad s=1,\dots,S_k  \, , \nonumber
\end{gather}
where $m$ and $g$ subscripts have been emitted for brevity. $I_{k,s+1}$ is the radiation intensity at the outgoing face of the $s^\text{th}$ segment of the $k^\text{th}$ characteristic, $\tilde{\varkappa}_{i}$ and ${Q}_{i}$ are the modified opacity and source in the $i^\text{th}$ material cell that is traced over by the $s^\text{th}$ segment. ${\bar{I}}_{k,s}$ is the segment-average radiation intensity, the optical thickness of each segment is
\begin{equation}
	\tau_{k,s} = \frac{\tilde{\varkappa}_{i}\ell_{k,s}}{\sin(\zeta)},\quad \alpha_{k,s} = \frac{1}{\tau_{k,s}} - \frac{e^{-\tau_{k,s}}}{1-e^{-\tau_{k,s}}} \, .
\end{equation}
The area of a characteristic segment is
\begin{equation}
	A_{k,s} = \ell_{k,s}w_{k} \, ,
\end{equation}
where $w_k$ is the width of the $k^\text{th}$ characteristic.
Radiation intensities on each characteristic are averaged with the segment widths and areas to find the cell-face and cell-average values, respectively, for $I$ on the underlying material spatial grid. These quantities are used to calculate
cell-  and face-average angular moments of high-order transport solution necessary to compute grid functions of the  Eddington tensor in material spatial cells. The cell-average Eddington tensor on the $i^\text{th}$ material grid cell is defined by
\begin{equation}
	\boldsymbol{\mathfrak{f}}_i = \frac{ \sum_{m=1}^{M}\sum_{k,s\in C(i,m)} \omega_mA_{k,s}(\boldsymbol{\Omega}_m\otimes\boldsymbol{\Omega}_m)I_{k,s} }{ \sum_{m=1}^{M}\sum_{k,s\in C(i,m)} \omega_mA_{k,s}I_{k,s} } \, ,
\end{equation}
where $\omega_m$ are the angular quadrature weights and $C(i,m)$ is the set of characteristic segments for the direction $\boldsymbol{\Omega}_m$ which trace over the $i^\text{th}$ material grid cell. The face-average Eddington tensor on the $f^\text{th}$ cell face in the material grid is
\begin{equation}
	\boldsymbol{\mathfrak{f}}_f = \frac{ \sum_{m=1}^{M}\sum_{k,s\in F(f,m)} \omega_mw_{k}(\boldsymbol{\Omega}_m\otimes\boldsymbol{\Omega}_m)I_{k,s} }{ \sum_{m=1}^{M}\sum_{k,s\in F(f,m)} \omega_mw_{k}I_{k,s} } \, ,
\end{equation}
where $F(f,m)$ is the set of characteristic segments for the direction $\boldsymbol{\Omega}_m$ whose upwind face intersects with the $f^\text{th}$ material grid cell face.

The characteristic grid is constructed as follows. For each direction, an initial set of characteristics are calculated such that no characteristic traces over a vertex in the material spatial grid. The result is a `non-uniform' grid of characteristics for each direction with varying widths $w$. We then enforce a maximum width $h_\text{moc}$. Any characteristics calculated from the initial ray tracing procedure that have a width $w>h_\text{moc}$ are split into 2 characteristics of width ${w}/{2}$. Characteristics are continuously split into even halves by width until none are wider than $h_\text{moc}$. Finally, each characteristic is formed into segments length-wise whose boundaries are at the intersection point of a given characteristic with a cell boundary in the material grid.

With this formulation, the TRT problem will be discretized with two unique computational grids in space.
This yields more flexibility in that both the underlying material grid and the high-order characteristic grid can be refined independently of one another.
If the characteristic grid is continuously refined on a static material grid, the closures provided by the BTE solution on the material grid are expected to converge to some value.
However the accuracy of these closures (and the low-order solution) will still be limited by the material grid.
When allocating computational resources in a simulation, it may be unclear which grid will benefit the most from extra refinement.
We perform convergence studies of each grid independently in Section \ref{sec:res} to quantify these effects.


\section{Numerical Results} \label{sec:res}

Analysis is performed using the classical Fleck-Cummings (F-C) test problem  in 2D Cartesian geometry \cite{fleck-1971}. The test domain is a homogeneous $6\times 6$ cm slab with material defined by the spectral opacity  $\varkappa_\nu = \frac{27}{\nu^3}(1-e^{-\nu/T})$ and energy density $\varepsilon=c_vT$ where $c_v=0.5917 a_R (T^\text{in})^3$. $T^\text{in}=1$ KeV is the temperature of radiation incoming to the domain on the left boundary, with the other boundaries being at vacuum. The slab is initially at a temperature of $T^0=1$ eV.
The F-C test is simulated for times $t\in[0,3\text{ns}]$. 150 uniform time steps are used $\Delta t = 2\times 10^{-2}$ ns. 17 frequency groups are used and 144 discrete directions. The Abu-Shumays angular quadrature set q461214 with
36 discrete directions per quadrant is used \cite{abu-shumays-2001}. All low-order equations are discretized with the backward-Euler time integrator and a second-order finite volumes scheme. We define a uniform, orthogonal material spatial grid whose cells are width $h_\text{mat}$.
\begin{figure}[ht!]
	\centering
	\includegraphics[width=0.5\textwidth]{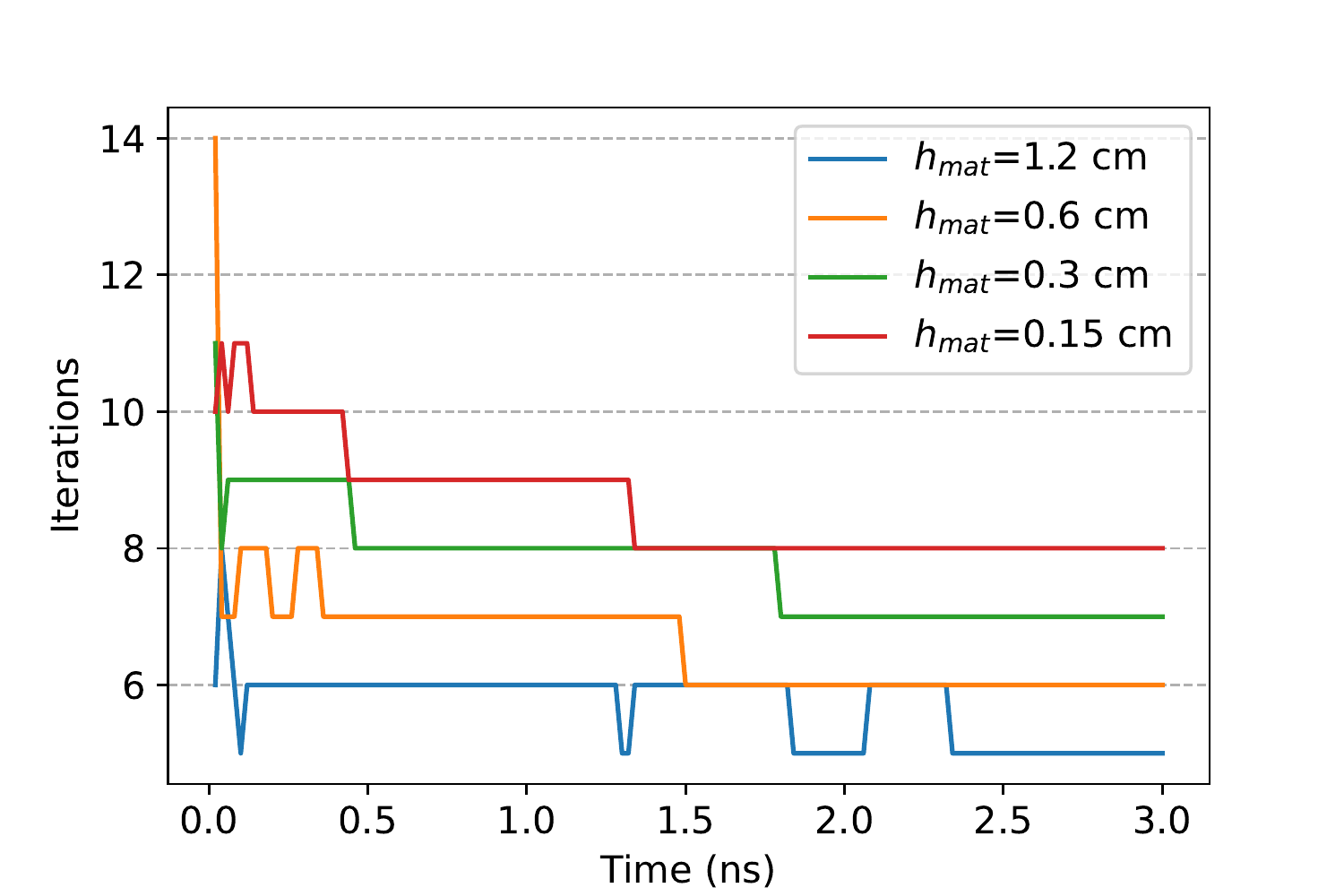}
	\caption{Iterations required with the MLQD algorithm to reach a convergence criteria of $\epsilon=10^{-14}$ at each time step while $h_\text{moc}=10^{-3}$ cm and $h_\text{mat}=1.2,0.6,0.3,0.015$ cm. \label{fig:iterations}}
\end{figure}

Figure \ref{fig:iterations} shows the number of high-order transport iterations required with the MLQD algorithm at each time step to reach a very tight convergence criteria of $\epsilon=10^{-14}$. Such a tight criteria is enforced to ensure no effects are missed in the iterative behavior,
and comparisons of the solution on different spatial grid sizes can be made to high precision. Results are shown for $h_\text{moc}=10^{-3}$ cm with several values for $h_\text{mat}$. Iterations increase slightly with refinements in $h_\text{mat}$. For all considered mesh widths the iterations converge rapidly. Here $h_\text{moc}$ is refined to a high level for all considered material grids. For $h_\text{mat}=0.15$ cm, there will be more than $150$ characteristics traced over every material grid cell for each discrete direction. In this way the characteristic grid can be considered static for all considered $h_\text{mat}$.
We note that iteration counts per time step change insignificantly with changes in $h_\text{moc}$.

Tables \ref{tab:matconv}, \ref{tab:mocconv} and \ref{tab:mocconv2} display the results of mesh refinement studies conducted for both $h_\text{mat}$ and $h_\text{moc}$. Each refinement study is conducted with one of the mesh widths held static. In Table \ref{tab:matconv}, $h_\text{moc}$ is held constant at $10^{-3}$ cm while $h_\text{mat}$ is refined in several steps. In Table \ref{tab:mocconv}, $h_\text{mat}$ is held constant at $0.15$ cm while $h_\text{moc}$ is refined in several steps. Table \ref{tab:mocconv2} holds $h_\text{mat}$ at 0.6 cm while $h_\text{moc}$ is refined.  Here $\|\Delta y_h\| = \| y_{h} - y_{2h} \|_{L_2}$ and $\rho_h^y$ = $\frac{\|\Delta y_{2h}\|}{\|\Delta y_{h}\|}$. $\rho_h^y$ signifies the estimated convergence rate of the variable $y$ with mesh refinement. For Table \ref{tab:matconv}, $h_\text{moc}=10^{-3}$ cm was chosen so that refinements in $h_\text{mat}$ would not cause implicit refinements in the characteristic grid. This allows for analysis on the effect of strictly refining the material grid with an effectively constant characteristic grid. $h_\text{mat}=0.15$ cm was chosen as the material grid width in Table \ref{tab:mocconv} to give a reasonably fine mesh on which to study characteristic grid refinements. Similarly, $h_\text{mat}=0.6$ was chosen as a suitably coarse material grid for Table \ref{tab:mocconv2}.

\begin{table}[ht!]
	\centering
	\caption{Convergence behavior for $T$ and $E$ with refinements in $h_\text{mat}$ while $h_\text{moc}$ held constant at $10^{-3}$ cm \label{tab:matconv}}
\medskip
	\begin{tabular}{|c||c|c||c|c|}
		\hline
		$h_\text{mat}$ & $\|\Delta T_{h_\text{mat}}\|$ & $\rho_{h_\text{mat}}^T$ & $\|\Delta E_{h_\text{mat}}\|$ & $\rho_{h_\text{mat}}^E$\\ \hline
		0.6 & 19.1 & -      &  2.13$\times10^{-1}$ & -\\ \hline
		0.3 & 9.70 & 1.96  &  1.10$\times10^{-1}$ & 1.94\\ \hline
		0.15 & 4.88 & 1.99 & 5.69$\times10^{-2}$ & 1.97\\ \hline
	\end{tabular}
\medskip
	\centering
	\caption{Convergence behavior for $T$ and $E$ with refinements in $h_\text{moc}$ while $h_\text{mat}$ held constant at $0.15$ cm\label{tab:mocconv}}
\medskip
	\begin{tabular}{|c|c||c|c||c|c|}
		\hline
		$\frac{h_\text{mat}}{h_\text{moc}}$  & $h_\text{moc}$ & $\|\Delta T_{h_\text{moc}}\|$ & $\rho_{h_\text{moc}}^T$ & $\|\Delta E_{h_\text{moc}}\|$ & $\rho_{h_\text{moc}}^E$\\ \hline
		8 & $2\times 10^{-2}$ & $8.01\times 10^{-5}$ & - & $8.71\times 10^{-7}$ & -\\ \hline
		16 & $1\times 10^{-2}$ & $3.47\times 10^{-5}$ & 2.31 & $3.08\times 10^{-7}$ & 2.82\\ \hline
		32 & $5\times 10^{-3}$ & $2.58\times 10^{-5}$ & 1.35 & $2.01\times 10^{-7}$ & 1.53\\ \hline
		64 & $2.5\times 10^{-3}$ & $7.57\times 10^{-6}$ & 3.41 & $7.51\times 10^{-8}$ & 2.68\\ \hline
	\end{tabular}
\medskip
	\centering
	\caption{Convergence behavior for $T$ and $E$ with refinements in $h_\text{moc}$ while $h_\text{mat}$ held constant at $0.6$ cm\label{tab:mocconv2}}
\medskip
	\begin{tabular}{|c|c||c|c||c|c|}
		\hline
		$\frac{h_\text{mat}}{h_\text{moc}}$  & $h_\text{moc}$ & $\|\Delta T_{h_\text{moc}}\|$ & $\rho_{h_\text{moc}}^T$ & $\|\Delta E_{h_\text{moc}}\|$ & $\rho_{h_\text{moc}}^E$\\ \hline
		8 & $7.5\times 10^{-2}$ & $5.07\times 10^{-3}$ & - & $5.45\times 10^{-5}$ & -\\ \hline
		16 & $3.75\times 10^{-2}$ & $1.96\times 10^{-3}$ &  2.59 & $2.59\times 10^{-5}$ & 2.10\\ \hline
		32 & $1.875\times 10^{-2}$ & $1.04\times 10^{-3}$ & 1.88 & $1.05\times 10^{-5}$ & 2.48\\ \hline
		64 & $9.375\times 10^{-3}$ & $2.93\times 10^{-4}$ & 3.55 & $3.04\times 10^{-6}$ & 3.44\\ \hline
		128 & $4.6875\times 10^{-3}$ & $1.56\times 10^{-4}$ & 1.88 & $1.04\times 10^{-6}$ & 2.93\\ \hline
		256 & $2.34375\times 10^{-3}$ & $4.35\times 10^{-5}$ & 3.59 & $2.99\times 10^{-7}$ & 3.48\\ \hline
	\end{tabular}
\end{table}

\begin{figure}[ht!]
	\centering
	\includegraphics[width=0.45\textwidth]{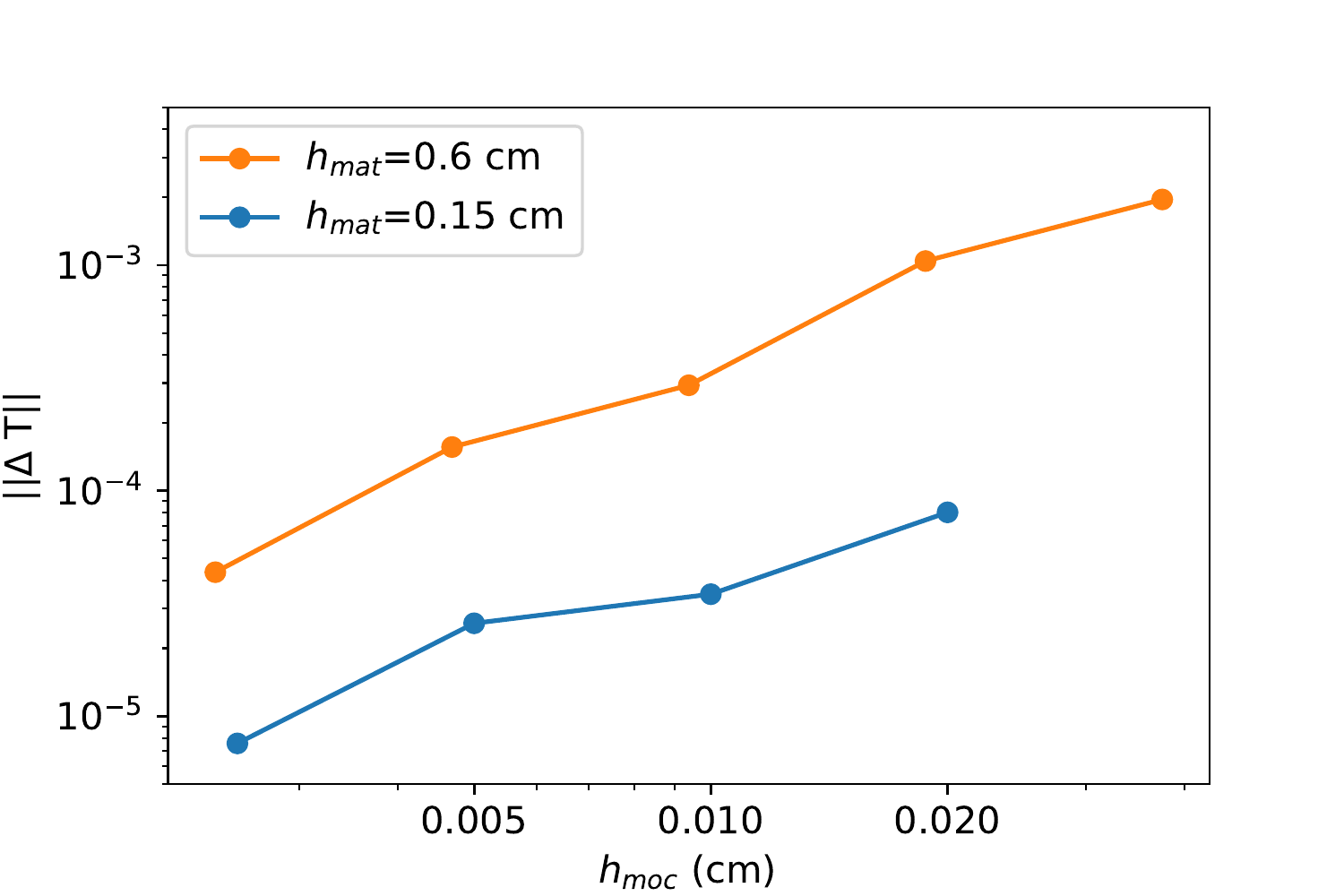}\\
	\includegraphics[width=0.45\textwidth]{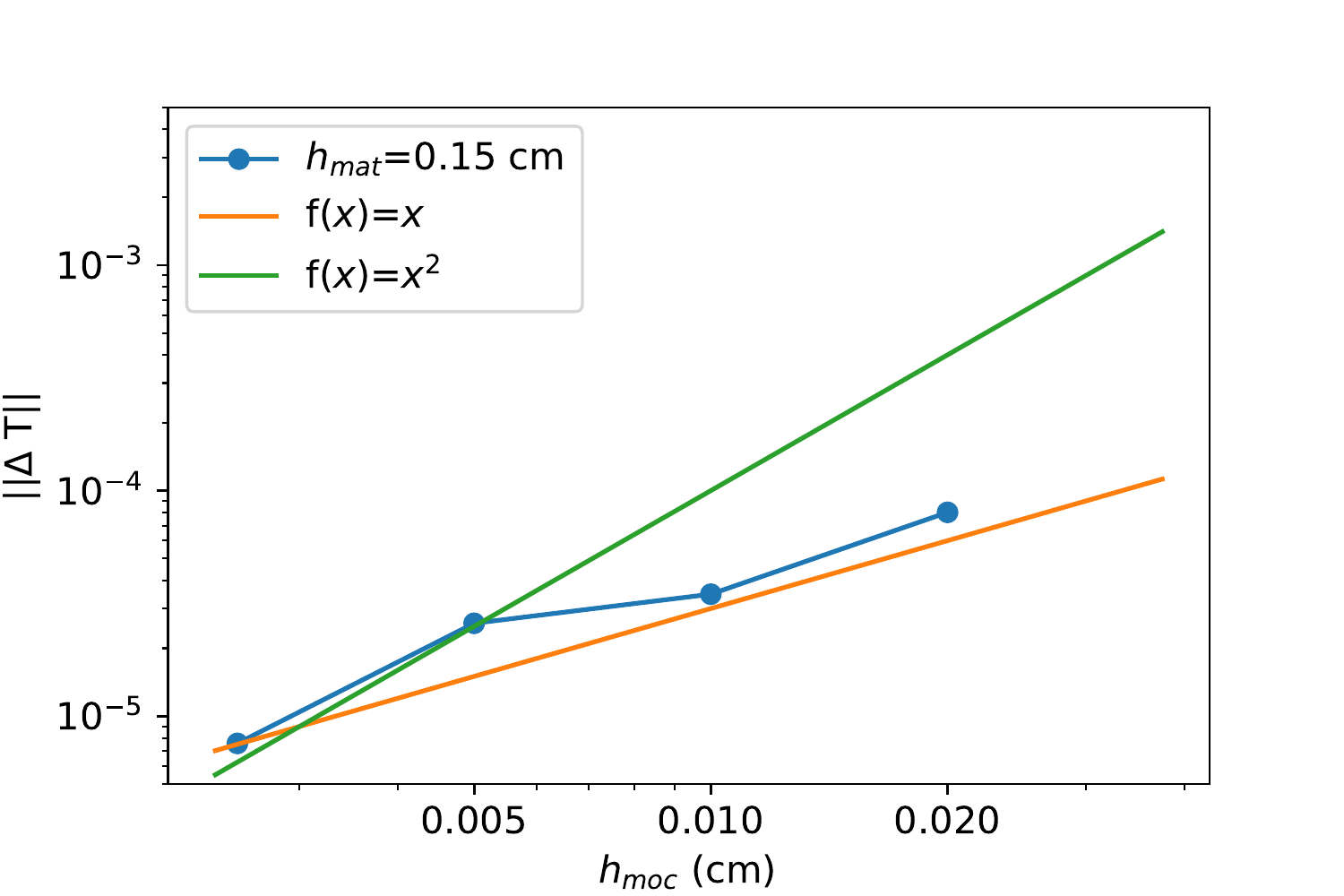}
	\includegraphics[width=0.45\textwidth]{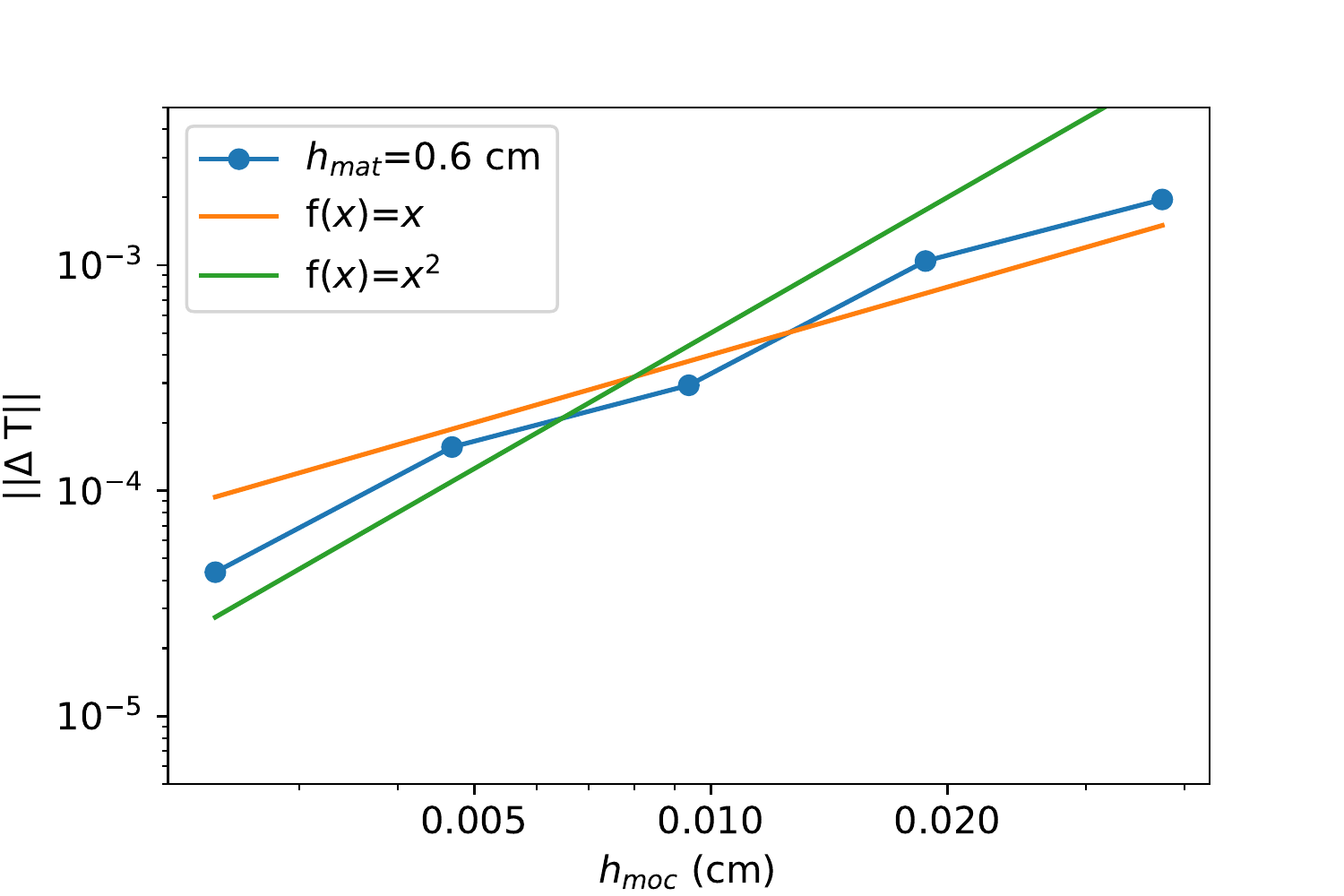}
	\caption{Plots of $\|\Delta T_{h_\text{moc}}\|$ vs $h_\text{moc}$ from Tables \ref{tab:mocconv} and \ref{tab:mocconv2} Top: $h_\text{mat}=0.6,\ 0.15$ cm; Bottom left: $h_\text{mat}=0.15$ with sample convergence rate lines; Bottom right: $h_\text{mat}=0.6$ with sample convergence rate lines. \label{fig:Tconv}}
\end{figure}

\begin{figure}[ht!]
	\centering
	\includegraphics[width=0.45\textwidth]{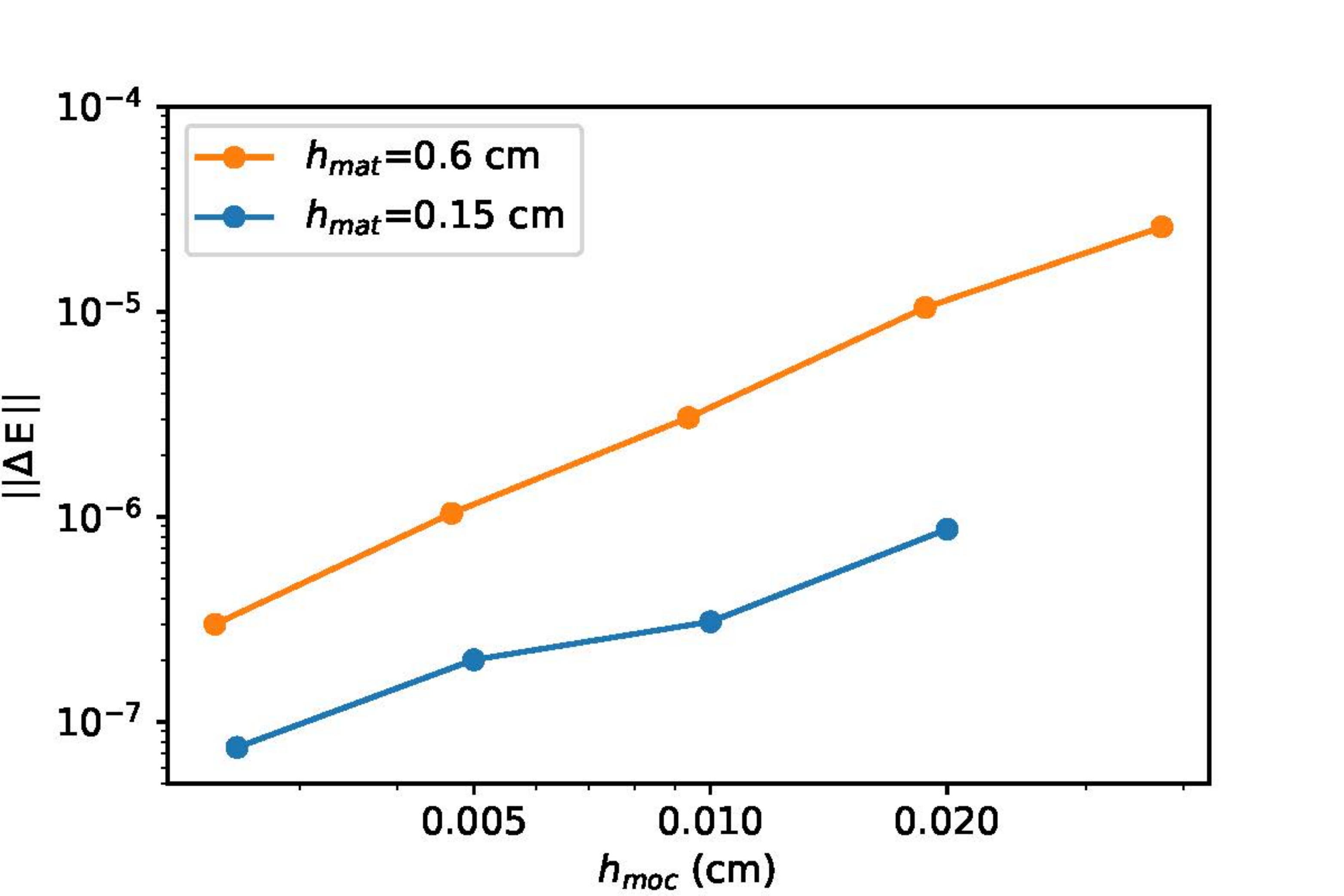}\\
	\includegraphics[width=0.45\textwidth]{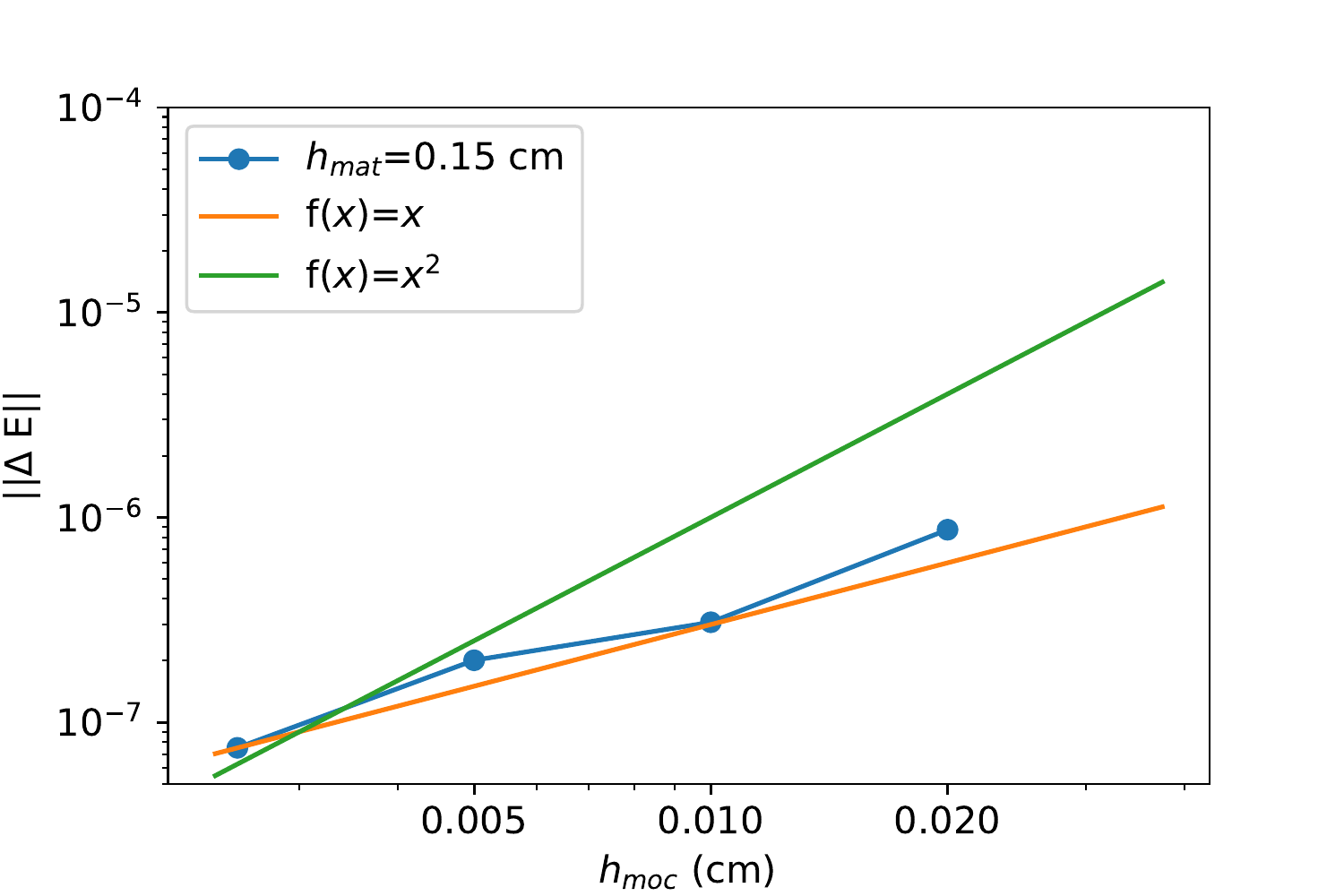}
	\includegraphics[width=0.45\textwidth]{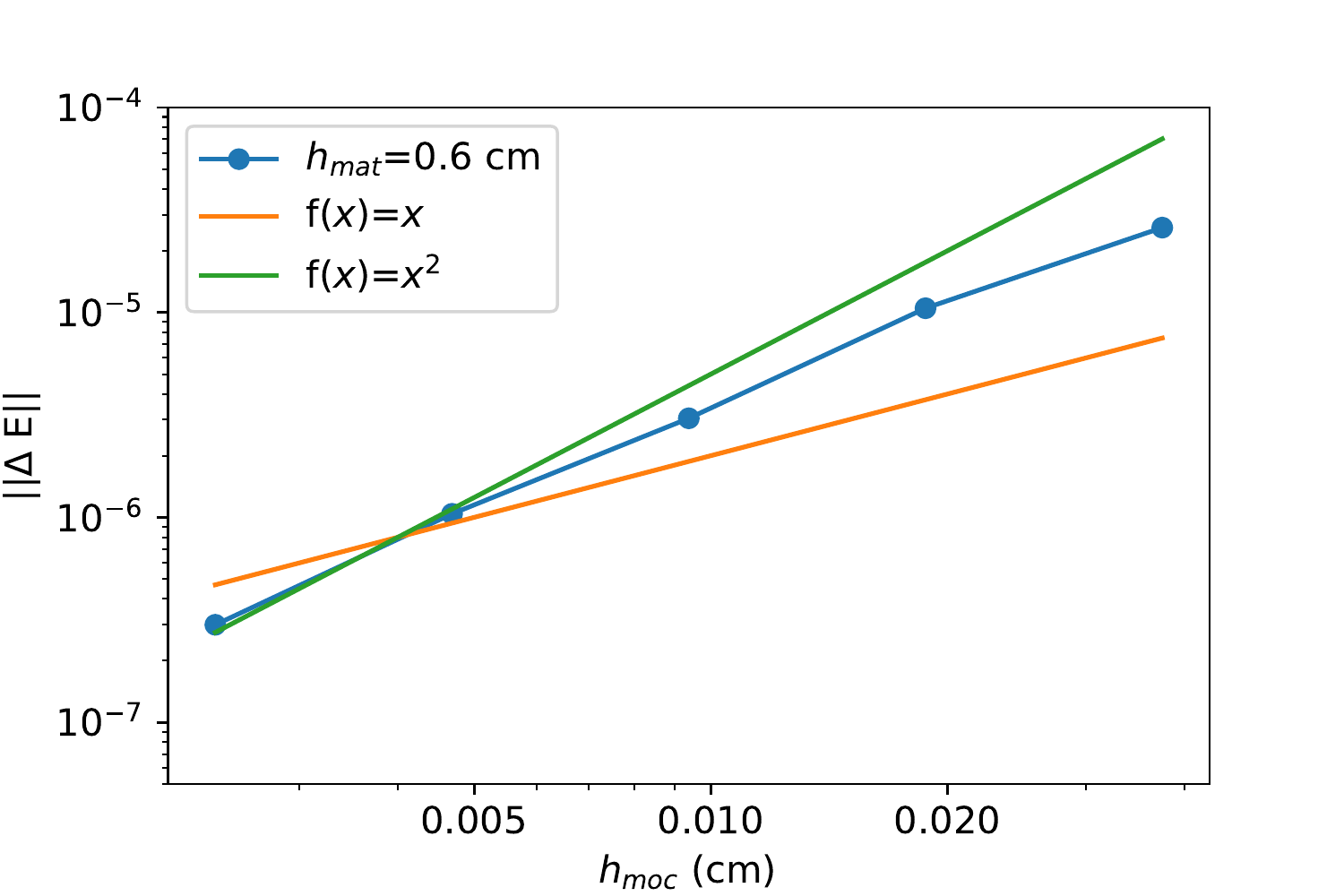}
	\caption{Plots of $\|\Delta E_{h_\text{moc}}\|$ vs $h_\text{moc}$ from Tables \ref{tab:mocconv} and \ref{tab:mocconv2} Top: $h_\text{mat}=0.6,\ 0.15$ cm; Bottom left: $h_\text{mat}=0.15$ with sample convergence rate lines; Bottom right: $h_\text{mat}=0.6$ with sample convergence rate lines. \label{fig:Econv}}
\end{figure}

The convergence rates for both $T$ and $E$ displayed for refinements in $h_\text{mat}$ in Table \ref{tab:matconv} are very close to first order. The convergence rates for refinements in $h_\text{moc}$ as shown in Tables \ref{tab:mocconv} and \ref{tab:mocconv2} are not as clear cut, taking on values between 1.3 and 3.6.
Figures \ref{fig:Tconv} and \ref{fig:Econv} display plots of $\|\Delta T_{h_\text{moc}}\|$ and $\|\Delta E_{h_\text{moc}}\|$ in log-log format. The values shown in both Tables \ref{tab:mocconv} and \ref{tab:mocconv2} are plotted together, and separately alongside sample lines representative of first and second order convergence rates.
Visually, the solution converges at a rate closer to first order when $h_\text{mat}=0.15$ cm, and closer to second order while $h_\text{mat}=0.6$ cm.
\begin{table}[ht!]
\centering
\caption{Relative differences in $T$ and $E$ from refinements in $h_\text{moc}$ while $h_\text{mat}$ held constant at $0.6$ cm\label{tab:mocconv2_rel}}
\medskip
\begin{tabular}{|c|c|c|c|}
	\hline
	$\frac{h_\text{mat}}{h_\text{moc}}$  & $h_\text{moc}$ & $\|\Delta T_{h_\text{moc}}\|_\text{rel}$ & $\|\Delta E_{h_\text{moc}}\|_\text{rel}$\\ \hline
	8 & $7.5\times 10^{-2}$ & $9.75\times 10^{-5}$ & $1.53\times 10^{-4}$\\ \hline
	16 & $3.75\times 10^{-2}$ & $3.76\times 10^{-5}$ & $7.31\times 10^{-5}$\\ \hline
	32 & $1.875\times 10^{-2}$ & $2.01\times 10^{-5}$ & $2.95\times 10^{-5}$\\ \hline
	64 & $9.375\times 10^{-3}$ & $5.64\times 10^{-6}$ & $8.57\times 10^{-6}$\\ \hline
	128 & $4.6875\times 10^{-3}$ & $3.00\times 10^{-6}$ & $2.93\times 10^{-6}$\\ \hline
	256 & $2.34375\times 10^{-3}$ & $8.36\times 10^{-7}$ & $8.43\times 10^{-7}$\\ \hline
\end{tabular}
\end{table}

The observed effects in $h_\text{moc}$ convergence behaviors
could partially be attributed to the fact that in constructing the characteristic mesh, only a maximum width of characteristics is enforced. Depending on how the initial grid is constructed, halving this maximum width criteria may not result in a strict doubling of the number of characteristics. Furthermore, this effect can change with the direction of motion $\boldsymbol{\Omega}$.
At the moment, there does not exist any rigorous theory which gives an expected rate of convergence for the solution with refinements in the characteristic mesh.
The results given here show that the scheme converges but we do not see asymptotic behavior for the rate of convergence. To fully understand the MOLC/RTS schemes for the class of problems at hand, more detailed analysis will be required.
\begin{table}[ht!]
\centering
\caption{Relative differences in $T$ and $E$ from refinements in $h_\text{moc}$ while $h_\text{mat}$ held constant at $0.15$ cm\label{tab:mocconv_rel}}
\medskip
\begin{tabular}{|c|c|c|c|}
	\hline
	$\frac{h_\text{mat}}{h_\text{moc}}$  & $h_\text{moc}$ & $\|\Delta T_{h_\text{moc}}\|_\text{rel}$ & $\|\Delta E_{h_\text{moc}}\|_\text{rel}$\\ \hline
	8 & $2\times 10^{-2}$ & $6.18\times 10^{-6}$ & $9.71\times 10^{-6}$\\ \hline
	16 & $1\times 10^{-2}$ & $2.68\times 10^{-6}$ & $3.43\times 10^{-6}$\\ \hline
	32 & $5\times 10^{-3}$ & $1.99\times 10^{-6}$ & $2.24\times 10^{-6}$\\ \hline
	64 & $2.5\times 10^{-3}$ & $5.84\times 10^{-7}$ & $8.37\times 10^{-7}$\\ \hline
\end{tabular}
\end{table}

Tables \ref{tab:mocconv2_rel} and \ref{tab:mocconv_rel} give the relative normed difference in solutions between subsequent $h_\text{moc}$ grids. The relative differences are calculated as $\|\Delta y_h\|_\text{rel} = {\| y_{h} - y_{2h} \|_{L_2}}\big/{\| y_{h} \|_{L_2}}$. The relative error norms in both $T$ and $E$ are shown to decrease by roughly an order of magnitude from $h_\text{mat}=0.6$ to $h_\text{mat}=0.15$ for the same ratio of $\frac{h_\text{mat}}{h_\text{moc}}$.
Instead if the two tables are compared at similar values for $h_\text{moc}$ widths, the relative difference values are similar to one another.
In any case, the magnitudes of both the absolute and relative error norms are small. Absolute values are capped at $5\times 10^{-3}$ and relative values are capped at $10^{-4}$.
In comparison to the results in Table \ref{tab:matconv}, the error norms are significantly smaller for refinements in the characteristic grid.

When considering how to formulate the computational mesh for a simulation, it is important to allocate resources such that the grid which limits solution accuracy is refined to the highest level. The results presented here demonstrate that, at least for the refinement levels examined, it is the underlying material grid which limits solution accuracy. Using $\frac{h_\text{mat}}{h_\text{moc}}=8$ for both $h_\text{mat}=0.6,\ 0.15$ cm is enough to ensure the characteristic grid is not the limiting grid in accuracy. It is reasonable to extrapolate that $\frac{h_\text{mat}}{h_\text{moc}}=8$ will continue to be satisfactory for finer values of $h_\text{mat}$ given that the error norms in Tables \ref{tab:mocconv} and \ref{tab:mocconv2} decrease with smaller $h_\text{mat}$ for the same $\frac{h_\text{mat}}{h_\text{moc}}$.


\section{Conclusions} \label{sec:conc}

In this paper the MLQD method with MOLC/RTS for the BTE is formulated for TRT problems and analyzed.
The method was tested on the 2D F-C test problem.
The analyzed method presents an independent discretization scheme for HEDP simulations which allows for straightforward refinement or coarsening of the discrete BTE solution separately from involved multiphysics equations.
As the BTE solution is refined, closures for the LOQD equations will become more resolved on the given material grid.
Iterations are shown to converge rapidly at every time step for several refinements of the material grid.
Independent mesh refinement studies were performed for the underlying material spatial grid and the characteristic grid.
First order convergence of the solution for $T$ and $E$ is observed for refinements in the material grid.
Refinements in the characteristic grid for the BTE do not give a clear convergence rate, instead alternating between sub- and super-linear rates for $T$ and $E$.
This could be due to the fact that refinements in the characteristic grid do not strictly double the number of characteristics.
A more detailed analysis of these effects should be partaken in the future to understand the convergence behavior of the MOLC/RTS in these problems.
The absolute difference between solutions for subsequent mesh refinements on the material grid was shown to be significantly larger than for refinements of the characteristic grid.

\section*{Acknowledgements}

Los Alamos Report LA-UR-22-32224.
This research project was funded by the Sandia National Laboratory, Light Speed Grand Challenge, LDRD, Strong Shock Thrust.
The work the first author (JMC)  was also supported by the U.S. Department of Energy through the Los Alamos National Laboratory. Los Alamos National Laboratory is operated by Triad National Security, LLC, for the National Nuclear Security Administration of U.S. Department of Energy (Contract No. 89233218CNA000001).
The content of the information does not necessarily reflect the position or the policy of the federal government, and no official endorsement should be inferred.

\bibliographystyle{elsarticle-num}
\bibliography{jmc-dya-MOLC-TRT-arXiv}

\begin{thebibliography}{10}
\expandafter\ifx\csname url\endcsname\relax
  \def\url#1{\texttt{#1}}\fi
\expandafter\ifx\csname urlprefix\endcsname\relax\def\urlprefix{URL }\fi
\expandafter\ifx\csname href\endcsname\relax
  \def\href#1#2{#2} \def\path#1{#1}\fi

\bibitem{zel-1966}
{Y. B. Zeldovich}, {Y. P. Razier}, Physics of Shock Waves and High Temperature
  Hydrodynamic Phenomena, Academic, New York, 1966.

\bibitem{drake-2006}
R.~P. Drake, High Energy Density Physics: Fundamentals, Inertial Fusion and
  Experimental Astrophysics, Springer, 2006.

\bibitem{Vladimirov-1958}
V.~Vladimirov, Numerical solution of the kinetic equation for a sphere,
  Computational Mathematics 3 (1958) 3, (in Russian).

\bibitem{goldin-1960}
{V. Ya. Gol'din}, Characteristic difference scheme for non-stationary kinetic
  equation, Soviet Mathematics Doklady 1 (1960) 902--906.

\bibitem{takeuchi-1969}
{K. Takeuchi}, A numerical method for solving the neutron transport equation in
  finite cylindrical geometry, Journal of Nuclear Science and Technology 6
  (1969) 446--473.

\bibitem{askew-1972}
{J. R. Askew}, A characteristic formulation of the neutron transport equation
  in complicated geometries, Tech. Rep. M1108, Atomic Energy Establishment
  (1972).

\bibitem{brough-chudley-1980}
{H. D. Brough}, {C. T. Chudley}, Characteristic ray solutions of the transport
  equation, in: {J. Lewis}, {M. Becker} (Eds.), Advances in Nuclear Science and
  Technology, Vol.~12, Plenum Press, 1980, pp. 1--30.

\bibitem{Alcouffe-1981}
{R. E.Alcouffe}, {E.W. Larsen}, A review of characteristic methods used to
  solve the linear transport equation, in: Proc. Int. Topl. Mtg. Advances in
  Mathematical Methods for the Solution of Nuclear Engineering Problems,
  Vol.~1, Munich, Germany, 1981, p.~3.

\bibitem{askew-roth-1982}
{J. R. Askew}, {M. J. Roth}, {WIMS-E}: {A} scheme for neutronics calculations,
  Tech. Rep. AEEW-R-1315, United Kingdom Atomic Energy Authority (1982).

\bibitem{casmo-1993}
{M. Edenius}, {K. Ekberg}, {B. H. Forss\'en}, {D. Knott}, {CASM0-4}, A Fuel
  Assembly Burnup Program, User's manual, {Studsvik/SOA-93/1}, {S}tudsvik of
  {A}merica, 1993.

\bibitem{kord-2002}
{K. S. Smith}, {J. D. Rhodes, III}, Full-core, {2-D}, {LWR} core calculations
  with {CASMO-4E}, in: Proceedings of PHYSOR 2002, Seoul, Korea, 2002, p. 12
  pp.

\bibitem{Sanchez-2000}
{R. Sanchez}, {A. Chetaine}, A synthetic acceleration for a two dimensional
  characteristic method in unstructured meshes, Nuclear Science and Engineering
  136 (2000) 122--139.

\bibitem{zika-adams-2000}
{M. Zika}, {M. Adams}, Acceleration for long-characteristics assembly-level
  transport problems, Nuclear Science and Engineering 134 (2000) 135--158.

\bibitem{Abel-2002}
{T. Abel}, {B. Wandelt}, Adaptive ray tracing for radiative transfer around
  point sources, Monthly Notices of the Royal Astronomical Society 330 (2002)
  L53--L56.

\bibitem{Pandya-2009}
{T. Pandya}, {M. Adams}, Method of long characteristics applied in space and
  time, in: Proc. of International Conference on Mathematics, Computational
  Methods \& Reactor Physics (M\&C 2009), Saratoga Springs, NY, 2009.

\bibitem{Buntemeyer-2016}
{L. Buntemeyer}, {R. Banerjeer}, {T. Peters r}, {M. Klassendr}, {R. Pudritz},
  Radiation hydrodynamics using characteristics on adaptive decomposed domains
  for massively parallel star formation simulations, New Astronomy 43 (2016)
  49--69.

\bibitem{gol'din-1964}
{V. Ya. Gol'din}, A quasi-diffusion method of solving the kinetic equation,
  USSR Comp.\ Math.\ and Math.\ Phys. 4 (1964) 136--149.

\bibitem{auer-mihalas-1970}
L.~H. Auer, D.~Mihalas, On the use of variable {E}ddington factors in non-{LTE}
  stellar atmospheres computations, Monthly Notices of the Royal Astronomical
  Society 149 (1970) 65--74.

\bibitem{gol'din-1972}
{V. Ya. Gol'din}, {B. N. Chetverushkin}, Methods of solving one-dimensional
  problems of radiation gas dynamics, USSR Comp.\ Math.\ and Math.\ Phys. 12
  (1972) 177--189.

\bibitem{winkler-norman-mihalas-85}
{K.-H. A. Winkler, M. L. Norman and D. Mihalas}, Implicit adaptive-grid
  radiation hydrodynamics, in: Multiple Time Scales, Academic Press, 1985, pp.
  145--184.

\bibitem{PASE-1986}
{V. Ya. Gol'din}, {D. A. Gol'dina}, {A. V. Kolpakov}, {A. V. Shilkov},
  Mathematical modeling of hydrodynamics processes with high-energy density
  radiation, Problems of Atomic Sci. \& Eng.: Methods and Codes for Numerical
  Solution of Math. Physics Problems 2 (1986) 59--88, in Russian.

\bibitem{mihalas-FRH-1984}
{D. Mihalas}, {B. Weibel-Mihalas}, Foundation of Radiation Hydrodynamics,
  Oxford University Press, 1984.

\bibitem{dya-aristova-vya-mm1996}
D.~Y. Anistratov, E.~N. Aristova, V.~Y. Gol'din, A nonlinear method for solving
  problems of radiation transfer in a physical system, Mathematical Modeling 8
  (1996) 3--28, in Russian.

\bibitem{aristova-vya-avk-m&c1999}
{E. N. Aristova}, {V. Ya. Gol'din}, {A. V. Kolpakov}, Multidimensional
  calculations of radiation transport by nonlinear quasi-diffusion method, in:
  Proc. of Int. Conf. on Math. and Comp., M\&C 1999, Madrid, Spain, 1999, pp.
  667--676.

\bibitem{dya-jcp-2019}
D.~Y. Anistratov, Stability analysis of a multilevel quasidiffusion method for
  thermal radiative transfer problems, Journal of Computational Physics 376
  (2019) 186--209.

\bibitem{dya-vyag-ttsp}
{D. Y. Anistratov}, {V. Ya. Gol'din}, Nonlinear methods for solving particle
  transport problems, Transport Theory and Statistical Physics 22 (1993)
  42--77.

\bibitem{pg-dya-jcp-2020}
{P. Ghassemi}, {D. Y. Anistratov}, Multilevel quasidiffusion method with
  mixed-order time discretization for multigroup thermal radiative transfer
  problems, Journal of Computational Physics 409 (2020) 109315.

\bibitem{fleck-1971}
{J. A. Fleck}, {J. D. Cummings}, An implicit {M}onte {C}arlo scheme for
  calculating time and frequency dependent nonlinear radiation transport, J. of
  Comp. Phys. 8 (1971) 313--342.

\bibitem{abu-shumays-2001}
{L. K. Abu-Shumays}, Angular quadratures for improved transport computations,
  Transport Theory \& Statistical Physics 30 (2001) 169--204.

\end{thebibliography}


\end{document}